\documentclass{article} 
\usepackage{latexsym,amssymb}
\usepackage{hyperref}
\usepackage{graphicx}
\usepackage[frenchb,english]{babel}
\usepackage[matrix,arrow]{xy}
\def \be{\begin{eqnarray*}}
\def \ee{\end{eqnarray*}}
\def \ben{\begin{enumerate}}
\def \een{\end{enumerate}}
\def \beit{\begin{itemize}}
\def \eeit{\end{itemize}}
\def \bui#1#2{\mathrel{\mathop{\kern 0pt#1}\limits^{#2}}}
\def \buil#1#2{\mathrel{\mathop{\kern 0pt#1}\limits_{#2}}}
\def \bfll{\begin{flushleft}}
\def \efll{\end{flushleft}}
\def \bflr{\begin{flushright}}
\def \eflr{\end{flushright}}

\def \findemo{\hfill$\square$}
\def \lra{\longrightarrow}
\def \lmt{\longmapsto}

\def \wih{\widehat}
\def \wit{\widetilde}
\def \wnabla{\wit{\nabla}}
\def \cdotM{\buil{\cdot}{M}}
\def \cdotN{\buil{\cdot}{N}}

\def \diag{\mathrm{diag}}

\def \SO{\mathrm{SO}}

\def \Spin{\mathrm{Spin}}
\def \SU{\mathrm{SU}}
\def \su{\mathfrak{su}}

\def \tr{\mathrm{tr}}

\def \la{\langle}
\def \ra{\rangle}

\newcommand{\rquot}[2]{\raisebox{0.5ex}{$#1$}\!/\!\raisebox{-0.5ex}{$#2$}}
\newcommand{\pa}[1]{\left(#1\right)}

\newtheorem{ethm}{Theorem}[section]

\newtheorem{elemme}[ethm]{Lemma}

\newtheorem{ecor}[ethm]{Corollary}
\newtheorem{prop}[ethm]{Proposition}

\def\SO{{\rm{SO}}} 
\def\SU{{\rm{SU}}} 
\def \U{{\rm{U}}}
\newcommand{\HH}{\mathbb{H}} 
 
\newcommand{\C}{\mathbb{C}} 
\newcommand{\Z}{\mathbb{Z}} 
\newcommand{\R}{\mathbb{R}}

\newcommand{\s}{\mathbb{S}}

\def\tr{{\rm {tr}}}

\newcommand{\lto}{\ensuremath{\longrightarrow}} 
 
\newcommand{\function}[5] 
{\begin{eqnarray*}\begin{array}{r@{}ccl} 
#1\;\colon\;  & #2 &\lto & #3 \\[.05cm] 
  & #4 &\longmapsto  & #5 
\end{array}\end{eqnarray*} 
}

\newcommand{\CP}{\mathbb{C}\mathrm{P}}

\title{The spectrum of the twisted Dirac operator on K\"ahler submanifolds of the complex projective space}
\author{Nicolas Ginoux\footnote{Fakult\"at f\"ur Mathematik,
Universit\"at Regensburg,
D-93040 Regensburg,
E-mail: \texttt{nicolas.ginoux@mathematik.uni-regensburg.de}}\, and Georges Habib\footnote{Lebanese University, Faculty of Sciences II, Department of Mathematics, P.O. Box 90656 Fanar-Matn, Lebanon, 
E-mail: \texttt{ghabib@ul.edu.lb}}}
\begin{document} 

\maketitle

\begin{abstract}
\noindent
We establish an upper estimate for the small eigenvalues of the twisted Dirac operator on K\"ahler submanifolds in K\"ahler manifolds carrying K\"ahlerian Killing spinors.
We then compute the spectrum of the twisted Dirac operator of the canonical embedding $\C P^d\rightarrow \C P^n$ in order to test the sharpness of the upper bounds.  
\end{abstract}

\section{Introduction}

One of the basic tools to get upper bounds for the eigenvalues of the twisted Dirac operator on spin submanifolds is the min-max principle.
The idea consists in computing in terms of geometric quantities the so-called Rayleigh-quotient applied to some test section coming from the ambient manifold.
In \cite{Bar98}, C. B\"ar established with the help of the min-max principle upper eigenvalue estimates for submanifolds in $\R^{n+1},\, \s^{n+1}$ and $\HH^{n+1}$, estimate which is sharp in the first two cases.
In the same spirit, the first-named author studied in his PhD thesis \cite{Ginthese} different situations where the ambient manifold admits natural test-spinors carrying geometric information.\\

\noindent In this paper, we consider a closed spin K\"ahler submanifold $M$ of a K\"ahler spin manifold $\wit{M}$ and derive upper bounds for the small eigenvalues of the corresponding twisted Dirac operator in case $\wit{M}$ carries so-called K\"ahlerian Killing spinors (see (\ref{eq:defKKS}) for a definition).
Interestingly enough, the upper bound turns out to depend only on the complex dimension of $M$ (Theorem \ref{t:majalphareel}).
Whether this estimate is sharp is a much more involved question.
A first approach consists in finding lower bounds for the spectrum and to compare them with the upper ones.
In Section \ref{s:kirch}, we prove a Kirchberg-type lower bound for the eigenvalues of any twisted Dirac operator on a closed K\"ahler manifold (Corollary \ref{c:estimKirch}). 
Here the curvature of the twisting bundle has to be involved.
Even for the canonical embedding $\C P^d\rightarrow \C P^n$, the presence of that normal curvature does not allow to state the equality between the lower bound and the upper one, see Proposition \ref{p:RE}.
The next approach consists in computing explicitly the spectrum of the twisted Dirac operator, at least for particular embeddings.
In Section \ref{s:CPd}, we determine the eigenvalues (with multiplicities) of the twisted Dirac operator of the canonical embedding $\C P^d\rightarrow \C P^n$, using earlier results by M. Ben Halima \cite{BenHalima08}.
We first remark that the spinor bundle of the normal bundle splits into a direct sum of powers of the tautological bundle (Corollary \ref{c:SigmaTperp}).
We deduce the spectrum of the twisted Dirac operator in Theorem \ref{t:compspec}, where we also include the multiplicities with the help of Weyl's character formula.
We conclude that, for $d<\frac{n+1}{2}$, the twisted Dirac operator admits $0$ as a lowest eigenvalue and $(n+1)(2d+1-n)$ for $d\geq \frac{n+1}{2}$ (see Proposition \ref{pro:loweigen}).
This implies that, for $d=1$, the upper estimate is optimal for $n=3,5,7$, however it is no more optimal for $n\geq 9$.\\

\noindent This work is partially based on and extends the first-named author's PhD thesis \cite[Ch. 4]{Ginthese}.\\

\noindent{\bf\small Acknowledgment.}{\small We thank the Max-Planck Institute for Mathematics in the Sciences and the University of Regensburg for their support.}

\section{Upper bounds for the submanifold Dirac operator of a K\"ahler submanifold}
\setcounter{equation}{0}

\noindent In this section, we prove {\it a priori} upper bounds for the smallest eigenvalues of some twisted Dirac operator on complex submanifolds in K\"ahler manifolds admitting so-called K\"ahlerian Killing spinors.\\

\noindent Let $M^{2d}$ be an immersed almost-complex submanifold in a K\"ahler manifold $(\wit{M}^{2n},g,J)$ (``almost-complex'' means that $J(TM)=TM$). Then for the induced metric and almost-complex structure the manifold $(M^{2d},g,J)$ is K\"ahler, in particular its immersion is minimal in $(\wit{M}^{2n},g,J)$.
We denote by $\wit{\Omega}$, $\Omega$ and $\Omega_N$ the K\"ahler form of $(\wit{M}^{2n},g,J)$, $(M^{2d},g,J)$ and of the normal bundle $NM\lra M$ of the immersion respectively (in our convention, $\Omega(X,Y)=g(J(X),Y)$ for all $X,Y$).\\
Assuming both $(M^{2d},g,J)$ and $(\wit{M}^{2n},g,J)$ to be spin, the bundle $NM$ carries an induced spin structure such that the restricted (complex) spinor bundle $\Sigma\wit{M}_{|_M}$ of $\wit{M}$ can be identified with $\Sigma M\otimes\Sigma N$, where $\Sigma M$ and $\Sigma N$ are the spinor bundles of $M$ and $NM$ respectively.
Denote by ``$\cdotM$'', ``$\cdotN$'' and ``$\cdot$''the Clifford multiplications of $M$, $NM$ and $\wit{M}$ respectively.
By a suitable choice of invariant Hermitian inner product $\la\cdot\,,\cdot\ra$ (with associated norm $|\cdot|$) on $\Sigma\wit{M}$ the identification above can be made unitary.
Moreover, it can be assumed to respect the following rules: given any $X\in TM$ and $\nu\in NM$, one has
\begin{equation}\label{eq:relCliffmult} 
\left|\begin{array}{ll}X\cdot\varphi&=\{X\cdotM\otimes(\mathrm{Id}_{\Sigma^+N}-\mathrm{Id}_{\Sigma^-N})\}\varphi\\ \nu\cdot\varphi&=(\mathrm{Id}\otimes\nu\cdotN)\varphi,\end{array}\right.
\end{equation}
for all $\varphi\in \Sigma\wit{M}_{|_M}=\Sigma M\otimes\Sigma N$. Here $\Sigma N=\Sigma^+N\oplus\Sigma^-N$ stands for the orthogonal and parallel splitting induced by the complex volume form, see e.g. \cite[Sec. 1.2.1]{Ginthese} or \cite[Sec. 2.1]{GinMor02}.
The following Gauss-type formula holds for the spinorial Levi-Civita connections $\wnabla$ and $\nabla:=\nabla^{\Sigma M\otimes\Sigma N}$ on $\Sigma\wit{M}$ and $\Sigma M\otimes\Sigma N$ respectively: for all $X\in TM$ and $\varphi\in\Gamma(\Sigma\wit{M}_{|_M})$,
\begin{equation}\label{eq:Gauss}\wnabla_X\varphi=\nabla_X\varphi+\frac{1}{2}\sum_{j=1}^{2d}e_j\cdot II(X,e_j)\cdot\varphi, \end{equation}
where $(e_j)_{1\leq j\leq 2d}$ is any local orthonormal basis of $TM$ and $II$ the second fundamental form of the immersion.\\

\noindent Recall that, for a complex constant $\alpha$, an $\alpha$-K\"ahlerian Killing spinor on a  K\"ahler spin manifold $(\wit{M}^{2n},g,J)$ is a pair $(\psi,\phi)$ of spinors satisfying, for all $X\in T\wit{M}$,
\begin{equation}\label{eq:defKKS}
\left|\begin{array}{ll}
\wnabla_X\psi&=-\alpha p_-(X)\cdot\phi\\
\wnabla_X\phi&=-\alpha p_+(X)\cdot\psi,
      \end{array}\right.
\end{equation}
where $p_\pm(X):=\frac{1}{2}(X\mp iJ(X))$.
The existence of a non-zero $\alpha$-K\"ahlerian Killing spinor on $(\wit{M}^{2n},g,J)$ imposes the metric to be Einstein with scalar curvature $\wit{S}=4n(n+1)\alpha^2$ (in particular $\alpha$ must be either real or purely imaginary), the complex dimension $n$ of $\wit{M}$ to be odd and the spinors $\psi,\phi$ to lie in particular eigenspaces of the Clifford action of $\wit{\Omega}$, namely 
\begin{equation}\label{eq:actionKaehlspinKKS}
\left|\begin{array}{ll}\wit{\Omega}\cdot\psi&=-i\psi\\
\wit{\Omega}\cdot\phi&=i\phi.\\
      \end{array}\right.
\end{equation}
Actually a K\"ahler spin manifold carries a non-zero $\alpha$-K\"ahlerian Killing spinor with $\alpha\in\R^\times$ if and only if it is the twistor-space of a quaternionic-K\"ahler manifold with positive scalar curvature (in particular it must be $\C\mathrm{P}^n$ if $n\equiv 1\;(4)$),  see \cite{Moroi95}. For purely imaginary $\alpha$ only partial results are known, the prominent examples being the complex hyperbolic space \cite[Thm. 13]{Kirch93} as well as doubly-warped products associated to some circle bundles over hyperk\"ahler manifolds \cite{GinHabibSem}.\\
We need the following lemma \cite[Lemme 4.4]{Ginthese}:

\begin{elemme}\label{l:D2psiphi}
Let $(M^{2d},g,J)$ be a K\"ahler spin submanifold of a K\"ahler spin manifold $(\wit{M}^{2n},g,J)$ and assume the existence of an $\alpha$-K\"ahlerian Killing spinor $(\psi,\phi)$ on $(\wit{M}^{2n},g,J)$.
Then 
\begin{equation}\label{eq:D2psiphi}
(D_M^{\Sigma N})^2(\psi+\phi)=(d+1)^2\alpha^2(\psi+\phi)+\alpha^2\Omega_N\cdot\Omega_N\cdot(\psi+\phi). 
\end{equation}
\end{elemme}

\noindent{\it Proof}: Fix a local orthonormal basis $(e_j)_{1\leq j\leq 2n}$ of $T\wit{M}_{|_M}$ with $e_j\in TM$ for all $1\leq j\leq 2d$ and $e_j\in NM$ for all $2d+1\leq j\leq 2n$. Introduce the auxiliary Dirac-type operator $\wih{D}:=\sum_{j=1}^{2d}e_j\cdot\wnabla_{e_j}:\Gamma(\Sigma\wit{M}_{|_M})\lra\Gamma(\Sigma\wit{M}_{|_M})$.
As a consequence of the Gauss-type formula (\ref{eq:Gauss}), the operators $\wih{D}^2$ and $(D_M^{\Sigma N})^2$ are related by \cite[Lemme 4.1]{Ginthese}
\[ \wih{D}^2\varphi=\pa{D_M^{\Sigma N}}^2\varphi-d^2|H|^2\varphi-d\sum_{j=1}^{2d}e_j\cdot\nabla_{e_j}^NH\cdot\varphi,\]
where $H:=\frac{1}{2d}\tr(II)$ is the mean curvature vector field of the immersion.
In particular $\wih{D}^2$ and $(D_M^{\Sigma N})^2$ coincide as soon as the mean curvature vector field of the immersion vanishes, condition which is fulfilled here.
Using $\sum_{j=1}^{2n}p_+(e_j)\cdot p_-(e_j)=i\wit{\Omega}-n$ and $\sum_{j=1}^{2n}p_-(e_j)\cdot p_+(e_j)=-i\wit{\Omega}-n$, we compute:
\be 
\wih{D}\psi&=&\sum_{j=1}^{2d}e_j\cdot\wnabla_{e_j}\psi\\
&\bui{=}{\rm (\ref{eq:defKKS})}&-\alpha\sum_{j=1}^{2d}e_j\cdot p_-(e_j)\cdot\phi\\
&=&-\alpha\sum_{j=1}^{2d}p_+(e_j)\cdot p_-(e_j)\cdot\phi\\
&=&-\alpha(i\Omega\cdot-d)\phi\\
&=&-\alpha(i\wit{\Omega}\cdot-d)\phi+i\alpha\Omega_N\cdot\phi\\
&\bui{=}{\rm (\ref{eq:actionKaehlspinKKS})}&(d+1)\alpha\phi+i\alpha\Omega_N\cdot\phi.
\ee
Similarly,
\be 
\wih{D}\phi&=&\sum_{j=1}^{2d}e_j\cdot\wnabla_{e_j}\phi\\
&\bui{=}{\rm (\ref{eq:defKKS})}&-\alpha\sum_{j=1}^{2d}e_j\cdot p_+(e_j)\cdot\psi\\
&=&-\alpha\sum_{j=1}^{2d}p_-(e_j)\cdot p_+(e_j)\cdot\psi\\
&=&-\alpha(-i\Omega\cdot-d)\psi\\
&=&\alpha(i\wit{\Omega}\cdot+d)\psi-i\alpha\Omega_N\cdot\psi\\
&\bui{=}{\rm (\ref{eq:actionKaehlspinKKS})}&(d+1)\alpha\psi-i\alpha\Omega_N\cdot\psi,
\ee
so that
\[\wih{D}(\psi+\phi)=(d+1)\alpha(\psi+\phi)+i\alpha\Omega_N\cdot(\phi-\psi).\]
To compute $\wih{D}^2(\psi+\phi)$ we need the commutator of $\Omega_N\cdot$ with $\wih{D}$.
For any $\varphi\in\Gamma(\Sigma\wit{M}_{|_M})$, one has
\be
\wih{D}(\Omega_N\cdot\varphi)&=&\sum_{j=1}^{2d}e_j\cdot\wnabla_{e_j}(\Omega_N\cdot\varphi)\\
&=&\sum_{j=1}^{2d}e_j\cdot\wnabla_{e_j}\Omega_N\cdot\varphi+e_j\cdot\Omega_N\cdot\wnabla_{e_j}\varphi\\
&=&\sum_{j=1}^{2d}\Omega_N\cdot e_j\cdot\wnabla_{e_j}\varphi+e_j\cdot\wnabla_{e_j}\Omega_N\cdot\varphi\\
&=&\Omega_N\cdot\wih{D}\varphi+\sum_{j=1}^{2d}e_j\cdot\wnabla_{e_j}\Omega_N\cdot\varphi,
\ee
with, for all $X,Y\in TM$ and $\nu\in NM$,
\be 
(\wnabla_X\Omega_N)(Y,\nu)&=&-\Omega_N(\wnabla_XY,\nu)\\
&=&-g(J(\wnabla_XY),\nu)\\
&=&-g(J(II(X,Y)),\nu),
\ee
so that
\be 
\sum_{j=1}^{2d}e_j\cdot\wnabla_{e_j}\Omega_N\cdot\varphi&=&-\sum_{j,k=1}^{2d}\sum_{l=2d+1}^{2n}g(J(II(e_j,e_k)),e_l)e_j\cdot e_k\cdot e_l\cdot\varphi\\
&=&-\sum_{j,k=1}^{2d}e_j\cdot e_k\cdot J(II(e_j,e_k))\cdot\varphi\\
&=&\sum_{j=1}^{2d}J(II(e_j,e_j))\cdot\varphi\\
&=&0,
\ee
since the immersion is minimal.
Hence $\wih{D}(\Omega_N\cdot\varphi)=\Omega_N\cdot\wih{D}\varphi$ and we deduce that
\be 
\wih{D}^2(\psi+\phi)&=&(d+1)\alpha\wih{D}(\psi+\phi)+i\alpha\wih{D}(\Omega_N\cdot(\phi-\psi))\\
&=&(d+1)^2\alpha^2(\psi+\phi)+i(d+1)\alpha^2\Omega_N\cdot(\phi-\psi)+i\alpha\Omega_N\cdot\wih{D}(\phi-\psi)\\
&=&(d+1)^2\alpha^2(\psi+\phi)+i(d+1)\alpha^2\Omega_N\cdot(\phi-\psi)+i\alpha\Omega_N\cdot((d+1)\alpha(\psi-\phi)-i\alpha\Omega_N\cdot(\psi+\phi))\\
&=&(d+1)^2\alpha^2(\psi+\phi)+\alpha^2\Omega_N\cdot\Omega_N\cdot(\psi+\phi),
\ee
which concludes the proof.
\findemo
$ $\\

\noindent Next we formulate the main theorem of this section.
Its proof requires some further notations.
Given any rank-$2k$-Hermitian spin bundle $E\lra M$ with metric connection preserving the complex structure, the Clifford action of the K\"ahler form $\Omega_E$ of $E$ splits the spinor bundle $\Sigma E$ of $E$ into the orthogonal and parallel sum
\begin{equation}\label{eq:splitSigmaE}
\Sigma E=\bigoplus_{r=0}^k\Sigma_rE,
\end{equation}
where $\Sigma_rE:=\mathrm{Ker}(\Omega_E\cdot-i(2r-k)\mathrm{Id})$ is a subbundle of complex rank $\left(\begin{array}{c}k\\ r\end{array}\right)$. Moreover, given any $V\in E$, one has $p_\pm(V)\cdot\Sigma_rE\subset\Sigma_{r\pm1}E$.

\begin{ethm}[see \protect{\cite[Thm. 4.2]{Ginthese}}]\label{t:majalphareel}
Let $(M^{2d},g,J)$ be a closed K\"ahler spin submanifold of a K\"ahler spin manifold $(\wit{M}^{2n},g,J)$ and consider the induced spin structure on the normal bundle.
Assume the existence of a complex $\mu$-dimensional space of non-zero $\alpha$-K\"ahlerian Killing spinor on $(\wit{M}^{2n},g,J)$ for some $\alpha\in\R^\times$.
Then there are $\mu$ eigenvalues $\lambda$ of $(D_M^{\Sigma N})^2$ satisfying 
\begin{equation}\label{eq:majalphareel}
\lambda\leq\left\{\begin{array}{ll}(d+1)^2\alpha^2&\textrm{ if }d\textrm{ is odd}\\\\ d(d+2)\alpha^2 &\textrm{ if }d\textrm{ is even.}\end{array}\right.
\end{equation}
If moreover {\rm (\ref{eq:majalphareel})} is an equality for the smallest eigenvalue $\lambda$ and some odd $d$, then $\sum_{j=1}^{2d}e_j\cdot II(X,e_j)\cdot\psi=\sum_{j=1}^{2d}e_j\cdot II(X,e_j)\cdot\phi=0$.
\end{ethm}

\noindent{\it Proof}: Let $(\psi,\phi)$ be a non-zero $\alpha$-K\"ahlerian Killing spinor on $(\wit{M}^{2n},g,J)$.
We evaluate the Rayleigh-quotient $\frac{\int_M\la(D_M^{\Sigma N})^2(\psi+\phi),\psi+\phi\ra v_g}{\int_M\la\psi+\phi,\psi+\phi\ra v_g}$ and apply the min-max principle.
It can be deduced from Lemma \ref{l:D2psiphi} that
\be
\la(D_M^{\Sigma N})^2(\psi+\phi),\psi+\phi\ra&=&(d+1)^2\alpha^2|\psi+\phi|^2+\alpha^2\la\Omega_N\cdot\Omega_N\cdot(\psi+\phi),\psi+\phi\ra\\
&=&(d+1)^2\alpha^2|\psi+\phi|^2-\alpha^2|\Omega_N\cdot(\psi+\phi)|^2.
\ee
Using (\ref{eq:splitSigmaE}) for $E=NM$ we observe that $|\Omega_N\cdot(\psi+\phi)|\geq |\psi+\phi|$ if $n-d$ is odd (i.e., if $d$ is even) and is nonnegative otherwise.
The inequality follows.\\
If $d$ is odd and (\ref{eq:majalphareel}) is an equality for the smallest eigenvalue, then $(D_M^{\Sigma N})^2(\psi+\phi)=(d+1)^2\alpha^2(\psi+\phi)$ and $\Omega_N\cdot(\psi+\phi)=0$. Since $\wit{\Omega}=\Omega\oplus\Omega_N$ one has $\Sigma_r\wit{M}_{|_M}=\bigoplus_{s=0}^r\Sigma_sM\otimes\Sigma_{r-s}M$ (where each component vanishes as soon as the index exceeds its allowed bounds), so that $\psi\in\Gamma(\Sigma_{\frac{d-1}{2}}M\otimes\Sigma_{\frac{n-d}{2}}N)$ and $\phi\in\Gamma(\Sigma_{\frac{d+1}{2}}M\otimes\Sigma_{\frac{n-d}{2}}N)$.
Coming back to the Gauss-type equation (\ref{eq:Gauss}), one obtains
\[\left|\begin{array}{ll}\nabla_X\psi&=-\alpha p_-(X)\cdot\phi-\frac{1}{2}\sum_{j=1}^{2d}e_j\cdot II(X,e_j)\cdot\psi\\\nabla_X\phi&=-\alpha p_+(X)\cdot\psi-\frac{1}{2}\sum_{j=1}^{2d}e_j\cdot II(X,e_j)\cdot\phi\end{array}\right.\]
for all $X\in TM$.
Looking more precisely at the components of each side of those identities, one notices that, pointwise, $\nabla_X\psi\in\Sigma_{\frac{d-1}{2}}M\otimes\Sigma_{\frac{n-d}{2}}N$ and, using (\ref{eq:relCliffmult}), that $p_-(X)\cdot\phi\in\Sigma_{\frac{d-1}{2}}M\otimes\Sigma_{\frac{n-d}{2}}N$.
But pointwise $\sum_{j=1}^{2d}e_j\cdot II(X,e_j)\cdot\psi\in(\Sigma_{\frac{d-3}{2}}M\otimes\Sigma_{\frac{n-d-2}{2}}N)\oplus(\Sigma_{\frac{d-3}{2}}M\otimes\Sigma_{\frac{n-d+2}{2}}N)\oplus(\Sigma_{\frac{d+1}{2}}M\otimes\Sigma_{\frac{n-d-2}{2}}N)\oplus(\Sigma_{\frac{d+1}{2}}M\otimes\Sigma_{\frac{n-d+2}{2}}N)$, in particular this term must vanish. Analogously one has $\sum_{j=1}^{2d}e_j\cdot II(X,e_j)\cdot\phi=0$. 
This concludes the proof.
\findemo
$ $\\

\noindent To test the sharpness of the estimate (\ref{eq:majalphareel}), we would like to first compare it to an {\it a priori} lower bound.
This is the object of the next section.

\section{Kirchberg-type lower bounds}\label{s:kirch}
\setcounter{equation}{0}

\noindent In this section, we aim at giving Kirchberg type estimates for any twisted Dirac operator on closed K\"ahler spin manifolds.
First consider a K\"ahler spin manifold $M$ of complex dimension $d$ and let $E$ be any rank $2k$-vector bundle over $M$ endowed with a metric connection. We def\mbox{}ine a connection on the vector bundle $\Sigma:=\Sigma M\otimes E$ by $\nabla:=\nabla^{\Sigma M\otimes E}$.
The Dirac operator of $M$ twisted with $E$ is defined by $D_M^E:\Gamma(\Sigma)\rightarrow \Gamma(\Sigma)$, $D_M^E:=\sum_{i=1}^{2d}e_i\cdot\nabla_{e_i}$, where $\{e_i\}_{1\leq i\leq 2d}$ is any local orthonormal basis of $TM$ and ``$\cdot$'' stands for the Clifford multiplication tensorized with the identity of $E$.
The square of the Dirac-type operator $D_M^E$ is related to the rough Laplacian via the following Schr\"odinger-Lichnerowicz formula \cite[Thm. II.8.17]{LM}
$$(D_M^E)^2=\nabla^*\nabla +\frac{1}{4}({\rm Scal}_M+R^E),$$ 
where ${\rm Scal}_M$ denotes the scalar curvature of $M$ and $R^E$ is the endomorphism tensor f\mbox{}ield given by 
\function{R^E}{\Sigma}{\Sigma}{\psi}{2\sum_{i,j=1}^{2d}(e_i\cdot e_j\cdot {\rm Id}\otimes R^E_{e_i,e_j})\psi.}  


\noindent Recall that for any eigenvalue $\lambda$ of the Dirac operator, there exists an eigenspinor $\varphi$ associated with $\lambda$ such that $\varphi=\varphi_r+\varphi_{r+1},$ where $\varphi_r$ is a section in $\Sigma_{r}:=\Sigma_r M\otimes E.$
Here $\Sigma_r M$ is the subundle $\mathrm{Ker}(\Omega\cdot-i(2r-d)\mathrm{Id})$ of $\Sigma M$.
Such an eigenspinor $\varphi$ is called of type $(r,r+1)$.
In order to estimate the eigenvalues of the twisted Dirac operator we define, as in the classical way, on each subbundle $\Sigma_{r}$ the {\it twisted twistor operator} for all $X\in \Gamma(TM),\, \psi_r\in \Sigma_r$ by \cite{BHMM}
$$P_{X}\psi_{r}:=\nabla_{X}\psi_{r}+a_r p_{-}(X)\cdot D_{+}\psi_{r}+b_r p_{+}(X)\cdot D_{-}\psi_{r},$$ 
where $a_r=\frac{1}{2(r+1)},\,\, b_r=\frac{1}{2(m-r+1)}$ and  $D_{\pm}\psi_{r}=\sum_{i=1}^{2d}p_{\pm}(e_i)\cdot\nabla_{e_i}\psi_{r}.$\\\\ 
We state the following lemma:
\begin{elemme}\label{lem:esti}
For any eigenspinor $\varphi$ of type $(r,r+1)$, we have the following inequalities
\begin{eqnarray}
\lambda^{2}\geq 
\left\{\begin{array}{ll}
\frac{1}{4(1-a_r)}\mathop{\rm inf}\limits_{M_{\varphi_{r}}}({\rm Scal}_M+R^E_{\varphi_{r}}),\\\\
\frac{1}{4(1-b_{r+1})}\mathop{\rm inf}\limits_{M_{\varphi_{r+1}}}({\rm Scal}_M+R^E_{\varphi_{r+1}}),
\end{array}\right.
\label{eq:4}
\end{eqnarray}
where $R^E_\phi:=\Re(R^E(\phi),\frac{\phi}{|\phi|^2})$ is def\mbox{}ined on the set $M_\phi=\{x\in M|\,\, \phi(x)\neq 0\}$ for all spinor $\phi\in \Sigma.$ 
\end{elemme}

\noindent{\it Proof}: Using the identity $\sum_{i=1}^{2d}e_i\cdot P_{e_i}\psi_{r}=0$, one can easily prove by a straightforward computation that for any spinor $\psi_r\in \Sigma_r$   
\begin{equation}\label{eq:1}
|P\psi_{r}|^2=|\nabla\psi_{r}|^2-a_r|D_{+}\psi_{r}|^2-b_r|D_{-}\psi_{r}|^2.
\end{equation} 
Applying Equation (\ref{eq:1}) to $\varphi_{r}$ and $\varphi_{r+1}$ respectively and integrating over $M$, we get with the use of the Schr\"odinger-Lichnerowicz formula that 
$$0\leq \int_M [\lambda^2(1-a_r)-\frac{1}{4}({\rm Scal}_M+R^E_{\varphi_{r}})]|\varphi_{r}|^2.$$ 
Also that,  
$$0\leq \int_M [\lambda^2(1-b_{r+1})-\frac{1}{4}({\rm Scal}_M+R^E_{\varphi_{r+1}})]|\varphi_{r+1}|^2,$$ 
from which the proof of the lemma follows. 
\findemo\\\\ 

\noindent One can get rid of the dependence of the eigenspinors $\varphi_r$ and $\varphi_{r+1}$ in the r.h.s. of (\ref{eq:4}):

\begin{ecor}\label{c:estimKirch} Let $\kappa_1$ be the smallest eigenvalue of the (pointwise) self-adjoint operator $R^E.$ Then
\begin{eqnarray*}
\lambda^{2}\geq 
\left\{\begin{array}{ll}
\frac{d+1}{4d}({\rm Scal}_0+\kappa_1) &\textrm{ if }d\textrm{ is odd}\\\\
\frac{d}{4(d-1)}({\rm Scal}_0+\kappa_1) &\textrm{ if }d\textrm{ is even},
\end{array}\right.
\end{eqnarray*}
where ${\rm Scal}_0$ denotes the infimum of the scalar curvature.
\end{ecor}

\noindent{\it Proof}: Let us choose the lowest integer $r\in \{0,1,\cdots,d\}$ such that $\varphi$ is of type $(r,r+1).$
The existence of anti-linear parallel maps on $\Sigma M$ commuting with the Clifford multiplication (see e.g. \cite[Lemma 1]{GinLag}) allows to impose that  
$r\leq\frac{d-1}{2}$ if $d$ is odd and $r\leq\frac{d-2}{2}$ if $d$ is even. 
This concludes the proof. 
\findemo\\\\

\noindent In the following, we formulate the estimates (\ref{eq:4}) for the situation where $M$ is a complex submanifold of the projective space $\C{\rm P}^n$ and $E$ is the spinor bundle of the normal bundle $NM$ of the immersion.
To do this, we will estimate $R^E_\phi$ for all spinor field $\phi\in \Sigma$ in terms of the second fundamental form of the immersion.

\begin{prop}\label{p:RE} Let $(M^{2d},g,J)$ be a K\"ahler spin submanifold of the projective space $\C{\rm P}^n$. For all spinor field $\phi\in \Sigma,$ the curvature is equal to 
\begin{equation}\label{eq:RE1}
R^E_\phi=-4 \Re(\Omega\cdot\Omega_N\cdot\phi,\frac{\phi}{|\phi|^2})-\sum_{i,j,p=1}^{2d}\Re(e_i\cdot e_j\cdot II(e_i,e_p)\cdot II(e_j,e_p)\cdot \phi,\frac{\phi}{|\phi|^2})+|II|^2.
\end{equation} 
where $\Omega$ is the K\"ahler form of $M.$ 
\end{prop}

\noindent{\it Proof}: First, recall that for all $X,Y\in \Gamma(TM)$ and $U,V$ sections in $NM$, the normal curvature is related to the one of $\C{\rm P}^n$ via the formula \cite[Thm. 1.1.72]{Be}
\begin{eqnarray}
(R^{NM}_{X,Y}U,V)&=&(R^{\mathbb{CP}^n}_{X,Y}U,V)-(B_{X}U,B_{Y}V)+(B_{Y}U,B_{X}V)\nonumber\\ 
&=& 2g(X,J(Y))g(J(U),V)-\sum_{p=1}^{2d}g(II(X,e_p),U)g(II(Y,e_p),V)\nonumber\\&&+\sum_{p=1}^{2d}g(II(Y,e_p),U)g(II(X,e_p),V),
\label{eq:normcur}
\end{eqnarray}
where $B_X: NM\rightarrow TM$ is the tensor f\mbox{}ield def\mbox{}ined by $g(B_X U,Y)=-g(II(X,Y),U)$ and $\{e_p\}_{1\leq p\leq 2d}$ is a local orthonormal basis of $TM$. Here we used the fact that the curvature of $\C{\rm P}^n$ is given for all $X,Y,Z \in T\C{\rm P}^n$ by  
$$R^{\mathbb{CP}^n}_{X,Y}Z=(X\wedge Y+JX\wedge JY+2g(X,JY)J)Z$$ 
with $(X\wedge Y)Z=g(Y,Z)X-g(X,Z)Y$.
Hence by (\ref{eq:normcur}), the normal spinorial curvature associated with any spinor field $\phi$ is then equal to 
\begin{eqnarray*}
R^E_{e_i,e_j}\phi&=&\frac{1}{4}\sum_{k,l=1}^{2(n-d)}g(R^{NM}_{e_i,e_j}e_k,e_l)e_k\cdot e_l\cdot \phi\\
&=&\frac{1}{2}\sum_{k=1}^{2(n-d)}g(e_i,J(e_j)) e_k\cdot Je_k\cdot\phi-\frac{1}{2}\sum_{p=1}^{2d}[II(e_i,e_p)\cdot II(e_j,e_p)\cdot+g(II(e_i,e_p),II(e_j,e_p))]\phi.
\end{eqnarray*}
Thus, we deduce
\begin{eqnarray*} 
R^E(\phi)&=&2\sum_{i,j=1}^{2d} J(e_j)\cdot e_j\cdot\Omega_N\cdot \phi-\sum_{i,j,p=1}^{2d}e_i\cdot e_j\cdot II(e_i,e_p)\cdot II(e_j,e_p)\cdot \phi\\&&-e_i\cdot e_j\cdot g(II(e_i,e_p),II(e_j,e_p))\phi\\ 
&=&-4 \Omega\cdot\Omega_N\cdot\phi-\sum_{i,j,p=1}^{2d}e_i\cdot e_j\cdot II(e_i,e_p)\cdot II(e_j,e_p)\cdot \phi+|II|^2\phi.
\end{eqnarray*}
Finally, the scalar product of the last equality with $\frac{\phi}{|\phi|^2}$ f\mbox{}inishes the proof. 
\findemo\\\\


\noindent
As we said in the proof of Corollary \ref{c:estimKirch}, the integer $r$ can be chosen such that $r\leq\frac{d-1}{2}$ if $d$ is odd and $r\leq\frac{d-2}{2}$ if $d$ is even.
However, we note that {\sl a priori} no such choice can be made for $s$ once $r$ has been fixed.
In particular, one cannot conclude that the smallest twisted Dirac eigenvalue of a totally geodesic $M$ in $\wit{M}$ is $(d+1)^2$, even in the ``simplest'' case where $M=\CP^d$ (the $d$-dimensional complex projective space).
To test the sharpness of the estimate (\ref{eq:majalphareel}), we compute in the following section the spectrum of $D_M^{\Sigma N}$ for $M=\CP^d$ canonically embedded in $\CP^n$.

\section{The spectrum of the twisted Dirac operator $D_M^{\Sigma N}$ on the complex projective space}\label{s:CPd}

\noindent In this section, we compute the spectrum of the Dirac operator of $\CP^d$ twisted with the spinor bundle of its normal bundle when considered as canonically embedded in $\CP^n$.
The eigenvalues will be deduced from M. Ben Halima's computations \cite[Thm. 1]{BenHalima08}.
We also need to compute the multiplicities in order to compare the upper bound in (\ref{eq:majalphareel}) with an eigenvalue which may be greater than the smallest one.
The results are gathered in Theorems \ref{t:specDgammadm} and \ref{t:compspec} below.

\subsection{The complex projective space as a symmetric space}
\setcounter{equation}{0}

\noindent Consider the $d$-dimensional complex projective space $\CP^d$ as the right quotient $\rquot{\SU_{d+1}}{{\rm S}(\U_d\times\U_1)}$, where ${\rm S}(\U_d\times\U_1):=\{\left(\begin{array}{ll}B&0\\0&\det(B)^{-1}\end{array}\right)\,|\,B\in\U_d\}$.
In this section we want to describe its tangent bundle and its normal bundle when canonically embedded into $\CP^n$ as homogeneous bundles, that is, as bundles associated to the ${\rm S}(\U_d\times\U_1)$-principal bundle $\SU_{d+1}\lra\CP^d$ via some linear representation of ${\rm S}(\U_d\times\U_1)$.
The one corresponding to the tangent bundle is called the \emph{isotropy representation} of the homogeneous space $\rquot{\SU_{d+1}}{{\rm S}(\U_d\times\U_1)}$.
To compute it explicitly we consider the following $\mathrm{Ad}({\rm S}(\U_d\times\U_1))$-invariant complementary subspace
\begin{equation}\label{eq:defm}
\mathfrak{m}:=\Big\{\left(\begin{array}{cccc}0 &\ldots &0& z_1\\ \vdots&{}&\vdots&\vdots\\0&\ldots&0& z_d\\-\bar{z_1}&\ldots& -\bar{z_d}&0\end{array}\right)\,|\,(z_1,\ldots,z_d)\in\C^d\Big\}
\end{equation}
to the Lie-Algebra $\mathfrak{h}$ of ${\rm S}(\U_d\times\U_1)$ in the Lie-algebra $\su_{d+1}=\{X\in\C(d+1)\,|\,X^*=-X\textrm{ and }\tr(X)=0\}$ and fix the (real) basis $(A_1,J(A_1),\ldots,A_d,J(A_d))$ of $\mathfrak{m}$, where:
\beit\item $(A_l)_{jk}=1$ if $(j,k)=(l,d+1)$, $-1$ if $(j,k)=(d+1,l)$ and $0$ otherwise;
\item $(J(A_l))_{jk}=i$ if $(j,k)=(l,d+1)$ or $(j,k)=(d+1,l)$ and $0$ otherwise.
\eeit
It is easy to check that $J$ defines a complex structure on $\mathfrak{m}$, which then makes $\mathfrak{m}$ into a $d$-dimensional complex vector space, and that $[\mathfrak{m},\mathfrak{m}]\subset\mathfrak{h}$.
In particular $\CP^d$ is a symmetric space.

\begin{elemme}\label{l:isotrepTCPd}
The isotropy representation of the symmetric space $\rquot{\SU_{d+1}}{{\rm S}(\U_d\times\U_1)}$ is given in the complex basis $(A_1,\ldots,A_d)$ of $\mathfrak{m}$ by:
\be 
\alpha:{\rm S}(\U_d\times\U_1)&\lra&\U_d\\
\left(\begin{array}{cc}B&0\\0&\det(B)^{-1}\end{array}\right)&\lmt&\det(B)\cdot B.
\ee
\end{elemme}

\noindent{\it Proof}: For $k\in\{1,\ldots,d\}$ and $B\in\U_d$ we compute
\be 
\mathrm{Ad}(\left(\begin{array}{cc}B&0\\0&\det(B)^{-1}\end{array}\right))(A_k)&=&\left(\begin{array}{cc}B&0\\0&\det(B)^{-1}\end{array}\right)\cdot A_k\cdot\left(\begin{array}{cc}B^*&0\\0&\det(B)\end{array}\right)\\
&=&\left(\begin{array}{cc}B&0\\0&\det(B)^{-1}\end{array}\right)\cdot\left(\begin{array}{cccc}0&\ldots&0&0\\\vdots& &\vdots&0\\0&\ldots&0&\det(B)\\\vdots& &\vdots&0\\-B_{k1}^*&\ldots&-B_{kd}^*&0\end{array}\right)\\
&=&\left(\begin{array}{cccc}0&\ldots&0&\det(B)B_{1k}\\\vdots& &\vdots&\vdots\\0&\ldots&0&\det(B)B_{dk}\\-\det(B)^{-1}B_{k1}^*&\ldots&-\det(B)^{-1}B_{kd}^*&0\end{array}\right)\\
&=&\sum_{j=1}^d\Re e(\det(B)B_{jk})A_j+\Im m(\det(B)B_{jk})J(A_j)\\
&=&\sum_{j=1}^d\det(B)B_{jk}A_j,
\ee
which gives the result.
\findemo
$ $\\

\noindent Recall that the tautological bundle of $\CP^d$ is the complex line bundle $\gamma_d\lra\CP^d$ defined by
\[\gamma_d:=\{([z],v)\,|\,[z]\in\CP^d\textrm{ and }v\in[z]\}.\]
It carries a canonical Hermitian metric defined by $\la([z],v),([z],v')\ra:=\la v,v'\ra$.

\begin{elemme}\label{l:isotrepTperpCPd}
The normal bundle $T^\perp\CP^d$ of the canonical embedding $\CP^d\rightarrow\CP^n$, $[z]\mapsto [z,0_{n-d}]$, is unitarily isomorphic to $\gamma_d^*\otimes\C^{n-d}$, where $\gamma_d\lra\CP^d$ is the tautological bundle of $\CP^d$ and $\C^{n-d}$ carries its canonical Hermitian inner product.
In particular, the homogeneous bundle $T^\perp\CP^d\rightarrow\CP^d$ is associated to the ${\rm S}(\U_d\times\U_1)$-principal bundle $\SU_{d+1}\lra\CP^d$ via the representation
\be 
\rho:{\rm S}(\U_d\times\U_1)&\lra&\U_{n-d}\\
\left(\begin{array}{cc}B&0\\0&\det(B)^{-1}\end{array}\right)&\lmt&\det(B)\mathrm{I}_{n-d}.
\ee
\end{elemme}

\noindent{\it Proof}: Consider the map
\be 
\CP^d\times\C^{n-d}&\bui{\lra}{\phi}&\gamma_d\otimes T^\perp\CP^d\\
([z],v)&\lmt&([z],z)\otimes d_z\pi(0_{d+1},v),
\ee
where $\pi:\C^{n+1}\lra\CP^n$ is the canonical projection.
It can be easily checked that $\phi$ is well-defined (the identity $\pi(\lambda z)=\pi(z)$ implies $d_z\pi=\lambda d_{\lambda z}\pi$) and is a unitary vector-bundle-isomorphism.
This shows the first statement.
Let $(e_1,\ldots,e_{d+1})$ denote the canonical basis of $\C^{d+1}$. 
The map
\be 
\SU_{d+1}\times\C&\lra&\gamma_d\\
(A,\lambda)&\lmt&([Ae_{d+1}],\lambda Ae_{d+1})
\ee
induces a complex vector-bundle-isomorphism $\rquot{\SU_{d+1}\times\C}{\rm S}(\U_d\times\U_1)\lra\gamma_d$, where the right action of ${\rm S}(\U_d\times\U_1)$ onto $\SU_{d+1}\times\C$ is given by $(A,\lambda)\cdot\left(\begin{array}{cc}B&0\\0&\det(B)^{-1}\end{array}\right):=(A\cdot\left(\begin{array}{cc}B&0\\0&\det(B)^{-1}\end{array}\right),\det(B)\lambda)$.
Thus $\gamma_d$ is isomorphic to the homogeneous bundle over $\CP^d$ which is associated to the ${\rm S}(\U_d\times\U_1)$-principal bundle $\SU_{d+1}\lra\CP^d$ via the representation ${\rm S(\U_d\times\U_1)}\rightarrow\U_1$, $\left(\begin{array}{cc}B&0\\0&\det(B)^{-1}\end{array}\right)\mapsto \det(B)^{-1}$.
This concludes the proof.
\findemo
$ $\\

\noindent Note in particular that $T^\perp\CP^d$ is not trivial (and hence not flat because of $\pi_1(\CP^d)=0$).

\subsection{Spin structures on $T\CP^d$ and $T^\perp\CP^d$}

\noindent From now on we assume that both $d$ and $n$ are odd integers.
Then both $T\CP^d$ and $T\CP^n$ are spin, in particular $T^\perp\CP^d$ is spin.
Since $\CP^d$ is simply-connected, there is a unique spin structure on $T\CP^d$ and on $T^\perp\CP^d$.
In this section we describe those spin structures as homogeneous spin structures.
For that purpose one looks for Lie-group-homomorphisms ${\rm S}(\U_d\times\U_1)\bui{\rightarrow}{\tilde{\alpha}}\Spin_{2d}$ and ${\rm S}(\U_d\times\U_1)\bui{\rightarrow}{\tilde{\rho}}\Spin_{2(n-d)}$ lifting $\alpha$ and $\rho$ through the non-trivial two-fold-covering map $\Spin_{2k}\bui{\lra}{\xi}\SO_{2k}$.\\
First we recall the existence for any positive integer $k$ of a Lie-group homomorphism $\U_k\bui{\lra}{j}\Spin_{2k}^c$ with $\xi^c\circ j=\iota$, where $\Spin_{2k}^c:=\rquot{\Spin_{2k}\times\U_1}{\Z_2}$ is the spin$^c$ group, $\xi^c:\Spin_{2k}^c\lra\SO_{2k}\times\U_1$, $[u,z]\mapsto(\xi(u),z^2)$ is the canonical two-fold-covering map and $\iota:\U_k\lra\SO_{2k}\times\U_1$, $A\mapsto(A_\R,\det(A))$.
The Lie-group homomorphism $j$ can be explicitly described on elements of $\U_k$ of diagonal form as:
\[j(\diag(e^{i\lambda_1},\ldots,e^{i\lambda_k}))=e^{\frac{i}{2}\left(\sum_{j=1}^k\lambda_j\right)}\cdot\wit{R}_{e_1,J(e_1)}(\frac{\lambda_1}{2})\cdot\ldots\cdot\wit{R}_{e_{k},J(e_{k})}(\frac{\lambda_k}{2}),\]
where $J$ is the canonical complex structure on $\C^k$ and, for any orthonormal system $\{v,w\}$ in $\R^{2k}$ and $\lambda\in\R$, the element $\wit{R}_{v,w}(\lambda)\in\Spin_{2k}$ is defined by
\[ \wit{R}_{v,w}(\lambda):=\cos(\lambda)+\sin(\lambda)v\cdot w.\] 
To keep the notations simple we denote by $j$ both such Lie-group-homomorphisms $\U_d\lra \Spin_{2d}^c$ and $\U_{n-d}\lra \Spin_{2(n-d)}^c$.

\begin{elemme}\label{l:homogspinstr}
Let $d<n$ be odd integers.
\ben\item The spin structure on $T\CP^d$ is associated to the ${\rm S}(\U_d\times\U_1)$-principal bundle $\SU_{d+1}\lra\CP^d$ via the Lie-group-homomorphism \be 
\tilde{\alpha}:{\rm S}(\U_d\times\U_1)&\lra&\Spin_{2d}\\
\left(\begin{array}{cc}B&0\\0&\det(B)^{-1}\end{array}\right)&\lmt&\det(B)^{-\frac{d+1}{2}}\cdot j\circ\alpha(\left(\begin{array}{cc}B&0\\0&\det(B)^{-1}\end{array}\right)).
\ee
\item The spin structure on $T^\perp\CP^d$ is associated to the ${\rm S}(\U_d\times\U_1)$-principal bundle $\SU_{d+1}\lra\CP^d$ via the Lie-group-homomorphism
\be 
\tilde{\rho}:{\rm S}(\U_d\times\U_1)&\lra&\Spin_{2(n-d)}\\
\left(\begin{array}{cc}B&0\\0&\det(B)^{-1}\end{array}\right)&\lmt&\det(B)^{-\frac{n-d}{2}}\cdot j\circ\rho(\left(\begin{array}{cc}B&0\\0&\det(B)^{-1}\end{array}\right)).
\ee
\een 
\end{elemme}

\noindent{\it Proof}: It suffices to prove the results for elements of ${\rm S}(\U_d\times\U_1)$ of diagonal form.
Indeed any e\-le\-ment of ${\rm S}(\U_d\times\U_1)$ is conjugated in $\SU_{d+1}$ to such a diagonal matrix.
Since $\SU_{d+1}$ is simply-connected the map $\SU_{d+1}\rightarrow\SO_{2k}\times\U_1$, $P\mapsto (PAP^{-1},\det(A))$ (where $A\in U_k$ is arbitrary), admits a lift through $\Spin_{2k}^c\bui{\lra}{\xi^c}\SO_{2k}\times\U_1$ which is uniquely determined by the image of one single point.
Therefore the lifts under consideration are uniquely determined on diagonal elements.\\
For $\theta_1,\ldots,\theta_d\in\R$ let $M_{\theta_1,\ldots,\theta_d}:=\diag(e^{i\theta_1},\ldots,e^{i\theta_d},e^{-i(\sum_{j=1}^d\theta_j)})\in{\rm S}(\U_d\times\U_1)$.
Then
\[u_{\theta_1,\ldots,\theta_d}:=\wit{R}_{e_1,J(e_1)}(\frac{\theta_1+\sum_{j=1}^d\theta_j}{2})\cdot\ldots\cdot\wit{R}_{e_{d},J(e_{d})}(\frac{\theta_d+\sum_{j=1}^d\theta_j}{2})\]
lies in $\Spin_{2d}$, only depends on $[\theta_1,\ldots,\theta_d]\in\rquot{\R^d}{2\pi\Z^d}$ (if some $\theta_k$ is replaced by $\theta_k+2m\pi$, then $u_{\theta_1,\ldots,\theta_d}$ is replaced by $(-1)^{m(d-1)}u_{\theta_1,\ldots,\theta_d}$, and $d-1$ is even) with $\xi(u_{\theta_1,\ldots,\theta_d})=\alpha(M_{\theta_1,\ldots,\theta_d})$.
Therefore $\tilde{\alpha}(M_{\theta_1,\ldots,\theta_d})=u_{\theta_1,\ldots,\theta_d}$.
Moreover, 
\be 
j\circ\alpha(M_{\theta_1,\ldots,\theta_d})&=&e^{\frac{i}{2}\left(\sum_{j=1}^d\theta_j+\sum_{k=1}^d\theta_k\right)}\cdot\wit{R}_{e_1,J(e_1)}(\frac{\theta_1+\sum_{j=1}^d\theta_j}{2})\cdot\ldots\cdot\wit{R}_{e_{d},J(e_{d})}(\frac{\theta_d+\sum_{j=1}^d\theta_j}{2})\\
&=&e^{\frac{i(d+1)}{2}\sum_{j=1}^d\theta_j}\cdot\tilde{\alpha}(M_{\theta_1,\ldots,\theta_d})\\
&=&\det(\diag(e^{i\theta_1},\ldots,e^{i\theta_d}))^{\frac{d+1}{2}}\cdot\tilde{\alpha}(M_{\theta_1,\ldots,\theta_d}),
\ee
which proves $1.$\\
The other case is much the same: setting
\[ \tilde{\rho}(M_{\theta_1,\ldots,\theta_d}):=\wit{R}_{e_1,J(e_1)}(\frac{\sum_{j=1}^d\theta_j}{2})\cdot\ldots\cdot\wit{R}_{e_{n-d},J(e_{n-d})}(\frac{\sum_{j=1}^d\theta_j}{2}),\]
one obtains a well-defined Lie-group-homomorphism ${\rm S}(\U_d\times\U_1)\bui{\rightarrow}{\tilde{\rho}}\Spin_{2(n-d)}$ with $\xi\circ\tilde{\rho}=\rho$ (the integer $n-d$ is even) and
\be 
j\circ\rho(M_{\theta_1,\ldots,\theta_d})&=&e^{\frac{i}{2}\sum_{j=1}^{n-d}\sum_{k=1}^d\theta_k}\cdot\wit{R}_{e_1,J(e_1)}(\frac{\sum_{j=1}^d\theta_j}{2})\cdot\ldots\cdot\wit{R}_{e_{n-d},J(e_{n-d})}(\frac{\sum_{j=1}^d\theta_j}{2})\\
&=&\det(\diag(e^{i\theta_1},\ldots,e^{i\theta_d}))^{\frac{n-d}{2}}\tilde{\rho}(M_{\theta_1,\ldots,\theta_d}),
\ee
which shows $2$ and concludes the proof.
\findemo
$ $\\

\noindent In particular, we obtain the following

\begin{ecor}\label{c:SigmaTperp}
Let $d<n$ be odd integers and consider the canonical embedding $\CP^d\rightarrow\CP^n$ as above.
Then there exists a unitary and parallel isomorphism
\[ \Sigma(T^\perp\CP^d)\cong\bigoplus_{s=0}^{n-d}\left(\begin{array}{c}n-d\\ s\end{array}\right)\cdot\gamma_d^{\frac{n-d}{2}-s},\]
where $\Sigma(T^\perp\CP^d)$ denotes the (complex) spinor bundle of $T^\perp\CP^d$ and, for each $s\in\{0,\ldots,n-d\}$, the factor $\left(\begin{array}{c}n-d\\ s\end{array}\right)$ stands for the multiplicity with which the line bundle $\gamma_d^{\frac{n-d}{2}-s}$ appears in the splitting.
\end{ecor}

\noindent{\it Proof}: By Lemma \ref{l:homogspinstr} and Lemma \ref{l:isotrepTperpCPd}, one has, for any $B\in \U_d$:
\be 
\tilde{\rho}(\left(\begin{array}{cc}B&0\\0&\det(B)^{-1}\end{array}\right))&=&\det(B)^{-\frac{n-d}{2}}\cdot j\circ\rho(\left(\begin{array}{cc}B&0\\0&\det(B)^{-1}\end{array}\right))\\
&=&\det(B)^{-\frac{n-d}{2}}\cdot j(\det(B)\mathrm{I}_{n-d}).
\ee
Now it is elementary to prove that, for any positive integer $k$, any $z\in\U_1$ and any $s\in\{0,\ldots,k\}$,
\[\delta_{2k}\circ j(z\cdot\mathrm{I}_k)_{|_{\Sigma_{2k}^{(s)}}}=z^s\cdot\mathrm{Id}_{\Sigma_{2k}^{(s)}},\]
where $\Sigma_{2k}^{(s)}$ is the eigenspace of the Clifford action of the K\"ahler form to the eigenvalue $i(2s-k)$ in the spinor space $\Sigma_{2k}$.
In particular $\Sigma_{2k}^{(s)}$ splits into the direct sum of $\mathrm{dim}_{\C}(\Sigma_{2k}^{(s)})$ copies of some one-dimensional representation, with $\mathrm{dim}_{\C}(\Sigma_{2k}^{(s)})=\left(\begin{array}{c} k\\ s\end{array}\right)$.
Since $\Sigma_{2k}=\oplus_{s=0}^k\Sigma_{2k}^{(s)}$, we obtain the following splitting:
\be 
\delta_{2(n-d)}\circ\tilde{\rho}&=&\bigoplus_{s=0}^{n-d}\det(\cdot)^{-(\frac{n-d}{2}-s)}\otimes\mathrm{Id}_{\Sigma_{2(n-d)}^{(s)}}\\
&=&\bigoplus_{s=0}^{n-d}\det(\cdot)^{-(\frac{n-d}{2}-s)}\otimes\mathbf{1}_{\C}^{\left(\begin{array}{c} n-d\\ s\end{array}\right)},
\ee
where $\det(\cdot):{\rm S}(\U_d\times\U_1)\rightarrow\U_1$, $\left(\begin{array}{cc}B&0\\0&\det(B)^{-1}\end{array}\right)\mapsto\det(B)$, the trivial representation on $\C$ is denoted by $\mathbf{1}_{\C}$ and ``$\mathbf{1}_{\C}^l$'' means that this representation appears with multiplicity $l$.
\findemo

\subsection{The twisted Dirac operator on $\CP^d$}

\noindent As a consequence of Corollary \ref{c:SigmaTperp}, the tensor product $\Sigma(T\CP^d)\otimes\Sigma(T^\perp\CP^d)$ splits into subbundles of the form $\Sigma(T\CP^d)\otimes\gamma_d^m$ for some integer $m$.
Since this splitting is orthogonal and parallel, it is also preserved by the corresponding twisted Dirac operator.
Hence it suffices to describe the Dirac operator of the twisted spinor bundle $\Sigma(T\CP^d)\otimes\gamma_d^m$ over $\CP^d$ as an infinite sum of matrices, where $m\in\mathbb{Z}$ is an arbitrary (non-necessarily positive) integer.
The Dirac eigenvalues of $\Sigma(T\CP^d)\otimes\gamma_d^m$ have been computed by M. Ben Halima in \cite[Thm. 1]{BenHalima08}. Indeed, we have

\begin{ethm}\label{l:valpDtordu}
For an odd integer $d$ let $\CP^d$ be endowed with its Fubini-Study metric of constant holomorphic sectional curvature $4$.
For an arbitrary $m\in\mathbb{Z}$ let the $m^{\textrm{th}}$ power $\gamma_d^m$ of the tautological bundle of $\CP^d$ be endowed with its canonical metric and connection.
Then the eigenvalues (without multiplicities) of the square of the Dirac operator of $\CP^d$ twisted by $\gamma_d^m$ are given by the following families:
\ben\item $2(r+l)\cdot(d+1+2(l-m-\epsilon))$, where $r\in\{1,\ldots,d-1\}$, $\epsilon\in\{0,1\}$ and $l\in\mathbb{N}$ with $l\geq\max(\epsilon,\frac{d+1}{2}-r+m)$.
\item $2l(2l+d-1-2m)$, where $l\in\mathbb{N}$, $l\geq\max(0,m+\frac{d+1}{2})$.
\item $2(d+l)(d+1+2(l-m))$, where $l\in\mathbb{N}$, $l\geq\max(0,m-\frac{d+1}{2})$.
\een
\end{ethm}
\,\,\,\,

\noindent The first family of eigenvalues corresponds to an irreducible representation of $\SU_{d+1}$ with highest weight given by \cite[Prop. 2]{BenHalima08}
\[ \hspace{-1cm}(r+2l-\frac{d-1}{2}-m-\epsilon,\underbrace{r+l-\frac{d-1}{2}-m,\ldots,r+l-\frac{d-1}{2}-m}_{r-1},r+l-\frac{d+1}{2}-m+\epsilon,\underbrace{r+l-\frac{d+1}{2}-m,\ldots,r+l-\frac{d+1}{2}-m}_{d-r-1}).\]
Similarly, the second family of eigenvalues corresponds to the highest weight 
\[(2l-\frac{d+1}{2}-m,\underbrace{l-\frac{d+1}{2}-m,\ldots,l-\frac{d+1}{2}-m}_{d-1}).\]
The last family of eigenvalues corresponds to 
\[(2l+\frac{d+1}{2}-m,\underbrace{l+\frac{d+1}{2}-m,\ldots,l+\frac{d+1}{2}-m}_{d-1}).\]

\noindent In the following, we will determine the multiplicities of the eigenvalues in Theorem \ref{l:valpDtordu}. Indeed, we have

\begin{elemme}\label{l:multvalp}
Let $d\geq1$ be an odd integer and $m\in\mathbb{Z}$.
\ben\item 
The multiplicities of the first family of the eigenvalues are equal to
\[\frac{d(\frac{d+1}{2}+r-m+2l-\epsilon)}{(r+l)(\frac{d+1}{2}-m+l-\epsilon)}\cdot\left(\begin{array}{c}d+l-\epsilon\\ d\end{array}\right)\cdot\left(\begin{array}{c}d-1\\ d-r-\epsilon\end{array}\right)\cdot\left(\begin{array}{c}\frac{d-1}{2}+r-m+l\\ d\end{array}\right).\]
\item 
For the second family, we have
\[
\prod_{k=2}^d(1+\frac{l}{k-1})\cdot(1+\frac{2l-\frac{d+1}{2}-m}{d})\cdot\prod_{j=2}^d(1+\frac{l-\frac{d+1}{2}-m}{d-j+1}).\]
\item 
For the last family of eigenvalues, the multiplicities are equal to
\[
\prod_{k=2}^d(1+\frac{l}{k-1})\cdot(1+\frac{2l+\frac{d+1}{2}-m}{d})\cdot\prod_{j=2}^d(1+\frac{l+\frac{d+1}{2}-m}{d-j+1}).\]
\een
\end{elemme}

\noindent In our convention, a product taken on an empty index-set is equal to $1$.\\

\noindent{\it Proof}: The required multiplicity can be computed with the help of the Weyl's character formula \cite{Baum94a}
\[ \prod_{\alpha\in\Delta_+}\Big(1+\frac{\la\lambda,\alpha\ra}{\la\delta_+,\alpha\ra}\Big),\]
where $\lambda$ is a highest weight of an irreducible $\SU_{d+1}$-representation and $\Delta_+$ is the set of positive roots, i.e. 
\[\Delta_+=\{\theta_j-\theta_k,\,\, 1\leq j<k\leq d,\,\, \theta_j+\sum_{k=1}^d\theta_k,\,\, 1\leq j\leq d\}\] 
 and $\delta_+=\sum_{k=1}^d(d-k+1)\theta_k$ is the half-sum of the positive roots of $\SU_{d+1}$, see \cite[p. 442]{BenHalima08}. Here the scalar product $<.,.>$ is the Riemannian metric on the dual of a maximal torus of $\SU_{d+1}$, which is defined by the following product of matrices $<\lambda,\lambda'>=\lambda.\beta.^t \lambda'$ where $\beta$ is the matrix given by $\frac{2}{d+1}\big(-1+(d+1)\delta_{jk}\big)_{1\leq j,k\leq d}.$ 
To compute the quotient in the Weyl's character formula, we treat the three cases separately:\\
1. Consider $\alpha$ of the form $\alpha=\theta_j-\theta_k$ for some $1\leq j<k\leq d$.
Note that this form for $\alpha$ can only exist if $d>1$.
We compute
\be 
\beta\cdot\alpha&=&\beta\cdot\left(\begin{array}{c}0\\\vdots\\0\\1\\0\\\vdots\\0\\-1\\0\\\vdots\\0\end{array}\right)\\
&=&\frac{2}{d+1}\left(\begin{array}{cccc}d&-1&\ldots&-1\\ -1&\ddots&&\vdots\\\vdots&&\ddots&-1\\-1&\ldots&-1&d\end{array}\right)\cdot\left(\begin{array}{c}0\\\vdots\\0\\1\\0\\\vdots\\0\\-1\\0\\\vdots\\0\end{array}\right)\\
&=&\frac{2}{d+1}\left(\begin{array}{c}0\\\vdots\\0\\d+1\\0\\\vdots\\0\\-d-1\\0\\\vdots\\0\end{array}\right)\\
&=&2(\theta_j-\theta_k).
\ee
Therefore,
\be
\la\delta_+,\alpha\ra&=&2(d,d-1,\ldots,1)\cdot\left(\begin{array}{c}0\\\vdots\\0\\1\\0\\\vdots\\0\\-1\\0\\\vdots\\0\end{array}\right)\\
&=&2(d-j+1-(d-k+1))\\
&=&2(k-j).
\ee
For the highest weight $\lambda$ corresponding to the first family of eigenvalues, we denote by $u'$ the first component, $\vec{u}_+$ the $r-1$ components, $u$ the $r$-components and by $\vec{u}_-$ the last $d-r-1$ components. Thus, we have
\be 
\la\lambda,\alpha\ra&=&2(u',\vec{u}_+,u,\vec{u}_-)\cdot\left(\begin{array}{c}0\\\vdots\\0\\1\\0\\\vdots\\0\\-1\\0\\\vdots\\0\end{array}\right)\\
&=&\left|\begin{array}{lll}2(u'-u_+)&\textrm{case}& j=1,\, k\in\{2,\ldots,r\}\\
2(u'-u)&\textrm{case}& j=1,\, k=r+1\\
2(u'-u_-)&\textrm{case}& j=1,\, k\in\{r+2,\ldots,d\}\\
0&\textrm{case}& j,k\in\{2,\ldots,r\}\\
2(u_+-u)&\textrm{case}& j\in\{2,\ldots,r\},\,k=r+1\\
2(u_+-u_-)&\textrm{case}& j\in\{2,\ldots,r\},\,k\in\{r+2,\ldots,d\}\\
2(u-u_-)&\textrm{case}& j=r+1,\,k\in\{r+2,\ldots,d\}\\
0&\textrm{case}& j,k\in\{r+2,\ldots,d\}\end{array}\right.\\
&=&\left|\begin{array}{lll}2(l-\epsilon)&\textrm{case}& j=1,\, k\in\{2,\ldots,r\}\\
2(l+1-2\epsilon)&\textrm{case}& j=1,\, k=r+1\\
2(l+1-\epsilon)&\textrm{case}& j=1,\, k\in\{r+2,\ldots,d\}\\
0&\textrm{case}& j,k\in\{2,\ldots,r\}\\
2(1-\epsilon)&\textrm{case}& j\in\{2,\ldots,r\},\,k=r+1\\
2&\textrm{case}& j\in\{2,\ldots,r\},\,k\in\{r+2,\ldots,d\}\\
2\epsilon&\textrm{case}& j=r+1,\,k\in\{r+2,\ldots,d\}\\
0&\textrm{case}& j,k\in\{r+2,\ldots,d\}.
\end{array}\right.
\ee
We obtain, for $\alpha=\theta_j-\theta_k$ with $1\leq j<k\leq d$:
\be 
1+\frac{\la\lambda,\alpha\ra}{\la\delta_+,\alpha\ra}&=&\left|\begin{array}{lll}\frac{l-\epsilon+k-j}{k-j}&\textrm{case}& j=1,\, k\in\{2,\ldots,r\}\\
\frac{l+1-2\epsilon+k-j}{k-j}&\textrm{case}& j=1,\, k=r+1\\
\frac{l+1-\epsilon+k-j}{k-j}&\textrm{case}& j=1,\, k\in\{r+2,\ldots,d\}\\
1&\textrm{case}& j,k\in\{2,\ldots,r\}\\
\frac{1-\epsilon+k-j}{k-j}&\textrm{case}& j\in\{2,\ldots,r\},\,k=r+1\\
\frac{1+k-j}{k-j}&\textrm{case}& j\in\{2,\ldots,r\},\,k\in\{r+2,\ldots,d\}\\
\frac{\epsilon+k-j}{k-j}&\textrm{case}& j=r+1,\,k\in\{r+2,\ldots,d\}\\
1&\textrm{case}& j,k\in\{r+2,\ldots,d\}.
\end{array}\right.
\ee
Now choose $\alpha=\theta_j+\sum_{k=1}^d\theta_k$ with $j\in\{1,\ldots,d\}$, then
\be 
\beta\cdot\alpha&=&\beta\cdot\left(\begin{array}{c}1\\\vdots\\1\\2\\1\\\vdots\\1\end{array}\right)\\
&=&\frac{2}{d+1}\left(\begin{array}{cccc}d&-1&\ldots&-1\\ -1&\ddots&&\vdots\\\vdots&&\ddots&-1\\-1&\ldots&-1&d\end{array}\right)\cdot\left(\begin{array}{c}1\\\vdots\\1\\2\\1\\\vdots\\1\end{array}\right)\\
&=&\frac{2}{d+1}\left(\begin{array}{c}0\\\vdots\\0\\2d-(d-1)\\0\\\vdots\\0\end{array}\right)\\
&=&2\theta_j.
\ee
Therefore,
\be
\la\delta_+,\alpha\ra&=&2(d,d-1,\ldots,1)\cdot\left(\begin{array}{c}0\\\vdots\\0\\1\\0\\\vdots\\0\end{array}\right)\\
&=&2(d-j+1).
\ee
Using the same notations as above, we compute 
\be 
\la\lambda,\alpha\ra&=&2(u',\vec{u}_+,u,\vec{u}_-)\cdot\left(\begin{array}{c}0\\\vdots\\0\\1\\0\\\vdots\\0\end{array}\right)\\
&=&\left|\begin{array}{lll}2u'&\textrm{case}& j=1\\
2u_+&\textrm{case}& j\in\{2,\ldots,r\}\\
2u&\textrm{case}& j=r+1\\
2u_-&\textrm{case}& j\in\{r+2,\ldots,d\}
\end{array}\right.\\
&=&\left|\begin{array}{lll}2(u_-+1+l-\epsilon)&\textrm{case}& j=1\\
2(u_-+1)&\textrm{case}& j\in\{2,\ldots,r\}\\
2(u_-+\epsilon)&\textrm{case}& j=r+1\\
2u_-&\textrm{case}& j\in\{r+2,\ldots,d\}.
\end{array}\right.
\ee
We obtain, for $\alpha=\theta_j+\sum_{k=1}^d\theta_k$ with $j\in\{1,\ldots,d\}$:
\be 
1+\frac{\la\lambda,\alpha\ra}{\la\delta_+,\alpha\ra}&=&\left|\begin{array}{lll}\frac{u_-+1+l-\epsilon+d-j+1}{d-j+1}&\textrm{case}& j=1\\
\frac{u_-+1+d-j+1}{d-j+1}&\textrm{case}& j\in\{2,\ldots,r\}\\
\frac{u_-+\epsilon+d-j+1}{d-j+1}&\textrm{case}& j=r+1\\
\frac{u_-+d-j+1}{d-j+1}&\textrm{case}& j\in\{r+2,\ldots,d\}.
\end{array}\right.
\ee
In order to compute the product we separate both cases $\epsilon=0$ and $\epsilon=1$.\\
$\bullet$ {\it Case $\epsilon=0$}: Then
\be 
\prod_{\alpha\in\Delta_+}\Big(1+\frac{\la\lambda,\alpha\ra}{\la\delta_+,\alpha\ra}\Big)&=&(\prod_{k=2}^r\frac{l+k-1}{k-1})\cdot(\prod_{k=r+1}^d\frac{l+k}{k-1})\cdot(\prod_{j=2}^r\prod_{k=r+1}^d\frac{k+1-j}{k-j})\cdot\\
& &\phantom{\prod_{k=2}^r\frac{l+k-1}{k-1}}\frac{u_-+l+d+1}{d}\cdot(\prod_{j=2}^r\frac{u_-+d-j+2}{d-j+1})\cdot(\prod_{j=r+1}^d\frac{u_-+d-j+1}{d-j+1})\\
&=&\frac{(l+1)\cdot\ldots\cdot(l+r-1)\cdot(l+r+1)\cdot\ldots\cdot(l+d)}{1\cdot2\cdot\ldots\cdot(d-1)}\cdot\\
& &\phantom{\prod_{k=2}^r\frac{l+k-1}{k-1}}(\prod_{j=2}^r\frac{(r+2-j)\cdot\ldots\cdot(d+1-j)}{(r+1-j)\cdot\ldots\cdot(d-j)})\cdot\frac{u_-+l+d+1}{d}\cdot\\
& &\phantom{\prod_{k=2}^r\frac{l+k-1}{k-1}}\frac{(u_-+d)\cdot\ldots\cdot(u_-+d-r+2)}{(d-1)\cdot\ldots\cdot(d-r+1)}\cdot\frac{(u_-+d-r)\cdot\ldots\cdot(u_-+1)}{(d-r)\cdot\ldots\cdot2\cdot1}\\
&=&\frac{d}{l+r}\cdot\frac{(l+d)!}{d!\cdot l!}\cdot(\prod_{j=2}^r\frac{d+1-j}{r+1-j})\cdot\frac{u_-+l+d+1}{u_-+d-r+1}\cdot\frac{(u_-+d)!}{d!\cdot u_-!}\\
&=&\frac{d}{l+r}\cdot\left(\begin{array}{c}l+d\\ d\end{array}\right)\cdot\frac{(d-1)\cdot\ldots\cdot(d+1-r)}{(r-1)\cdot\ldots\cdot2\cdot1}\cdot\frac{u_-+l+d+1}{u_-+d-r+1}\cdot\left(\begin{array}{c}u_-+d\\ d\end{array}\right)\\
&=&\frac{d(u_-+l+d+1)}{(l+r)(u_-+d-r+1)}\cdot\left(\begin{array}{c}l+d\\ d\end{array}\right)\cdot\left(\begin{array}{c}d-1\\ r-1\end{array}\right)\cdot\left(\begin{array}{c}u_-+d\\ d\end{array}\right),
\ee
which gives for the multiplicity in this case (replace $u_-$ by $r-\frac{d+1}{2}-m+l$):
\[\prod_{\alpha\in\Delta_+}\Big(1+\frac{\la\lambda,\alpha\ra}{\la\delta_+,\alpha\ra}\Big)=\frac{d(\frac{d+1}{2}+r-m+2l)}{(r+l)(\frac{d+1}{2}-m+l)}\cdot\left(\begin{array}{c}d+l\\ d\end{array}\right)\cdot\left(\begin{array}{c}d-1\\ d-r\end{array}\right)\cdot\left(\begin{array}{c}\frac{d-1}{2}+r-m+l\\ d\end{array}\right).\]\\

\noindent $\bullet$ {\it Case $\epsilon=1$}: Then
\be 
\prod_{\alpha\in\Delta_+}\Big(1+\frac{\la\lambda,\alpha\ra}{\la\delta_+,\alpha\ra}\Big)&=&(\prod_{k=2}^{r+1}\frac{l+k-2}{k-1})\cdot(\prod_{k=r+2}^d\frac{l+k-1}{k-1})\cdot(\prod_{j=2}^{r+1}\prod_{k=r+2}^d\frac{k+1-j}{k-j})\cdot\\
& &\phantom{\prod_{k=2}^r\frac{l+k-1}{k-1}}\frac{u_-+l+d}{d}\cdot(\prod_{j=2}^{r+1}\frac{u_-+d-j+2}{d-j+1})\cdot(\prod_{j=r+2}^d\frac{u_-+d-j+1}{d-j+1})\\
&=&\frac{l\cdot\ldots\cdot(l+r-1)\cdot(l+r+1)\cdot\ldots\cdot(l+d-1)}{1\cdot2\cdot\ldots\cdot(d-1)}\cdot\\
& &\phantom{\frac{l+k-1}{k-1}}(\prod_{j=2}^{r+1}\frac{(r+3-j)\cdot\ldots\cdot(d+1-j)}{(r+2-j)\cdot\ldots\cdot(d-j)})\cdot\frac{u_-+l+d}{d}\cdot\\
& &\phantom{\frac{l+k-1}{k-1}}\frac{(u_-+d)\cdot\ldots\cdot(u_-+d-r+1)}{(d-1)\cdot\ldots\cdot(d-r)}\cdot\frac{(u_-+d-r-1)\cdot\ldots\cdot(u_-+1)}{(d-r-1)\cdot\ldots\cdot2\cdot1}\\
&=&\frac{d}{l+r}\cdot\frac{(l+d-1)!}{d!\cdot(l-1)!}\cdot(\prod_{j=2}^{r+1}\frac{d+1-j}{r+2-j})\cdot \frac{u_-+l+d}{u_-+d-r}\cdot\frac{(u_-+d)!}{u_-!\cdot d!}\\
&=&\frac{d(u_-+l+d)}{(l+r)(u_-+d-r)}\cdot\frac{(l+d-1)!}{d!\cdot(l-1)!}\cdot\frac{(d-1)!}{r!\cdot(d-r-1)!}\cdot\frac{(u_-+d)!}{u_-!\cdot d!}\\
&=&\frac{d(u_-+l+d)}{(l+r)(u_-+d-r)}\cdot\left(\begin{array}{c}l+d-1\\ d\end{array}\right)\cdot\left(\begin{array}{c}d-1\\ r\end{array}\right)\cdot\left(\begin{array}{c}u_-+d\\ d\end{array}\right),
\ee
which, replacing $u_-$ by $r-\frac{d+1}{2}-m+l$, gives
\[\prod_{\alpha\in\Delta_+}\Big(1+\frac{\la\lambda,\alpha\ra}{\la\delta_+,\alpha\ra}\Big)= \frac{d(\frac{d-1}{2}+r-m+2l)}{(r+l)(\frac{d-1}{2}-m+l)}\cdot\left(\begin{array}{c}d+l-1\\ d\end{array}\right)\cdot\left(\begin{array}{c}d-1\\ d-r-1\end{array}\right)\cdot\left(\begin{array}{c}\frac{d-1}{2}+r-m+l\\ d\end{array}\right).\]
This shows 1.\\
2. Consider $\alpha$ of the form $\alpha=\theta_j-\theta_k$ for some $1\leq j<k\leq d$.
We have already shown in the first part that $\la\delta_+,\alpha\ra=2(k-j)$.
Let us denote by $v_+$ the first component and $v$ the $d-1$ components of the highest weight corresponding to the second family of eigenvalues. We have
\be 
\la\lambda,\alpha\ra&=&2(v_+,v,\ldots,v)\cdot\left(\begin{array}{c}0\\\vdots\\0\\1\\0\\\vdots\\0\\-1\\0\\\vdots\\0\end{array}\right)\\
&=&\left|\begin{array}{lll}2l&\textrm{case}& j=1\\0&\textrm{case}& j>1.\end{array}\right.
\ee
Choosing $\alpha=\theta_j+\sum_{k=1}^d\theta_k$ with $j\in\{1,\ldots,d\}$, we already know that 
$\la\delta_+,\alpha\ra=2(d-j+1)$.
Moreover,
\be 
\la\lambda,\alpha\ra&=&2(v_+,v,\ldots,v)\cdot\left(\begin{array}{c}0\\\vdots\\0\\1\\0\\\vdots\\0\end{array}\right)\\
&=&\left|\begin{array}{lll}2(v+l)&\textrm{case}& j=1\\2v&\textrm{case}& j>1.\end{array}\right.
\ee
Hence the product is given by
\[
\prod_{\alpha\in\Delta_+}\Big(1+\frac{\la\lambda,\alpha\ra}{\la\delta_+,\alpha\ra}\Big)=\prod_{k=2}^d(1+\frac{l}{k-1})\cdot(1+\frac{v+l}{d})\cdot\prod_{j=2}^d(1+\frac{v}{d-j+1}).\]
Of course only the central factor appears in case $d=1$.
Replacing $v$ by its respective value gives 2. and 3. and concludes the proof.
\findemo
$ $\\


\noindent As a consequence of Lemma \ref{l:valpDtordu} and Lemma \ref{l:multvalp}, we obtain the

\begin{ethm}\label{t:specDgammadm}
Let $d$ be a positive odd integer and $m\in\mathbb{Z}$ be arbitrary.
Denote by $\gamma_d$ the tautological bundle of $\CP^d$.
Then the spectrum of the square of the Dirac operator of $\CP^d$ twisted with $\gamma_d^m$ is given by the following family of eigenvalues:
\ben\item $2(r+l)\cdot(d+1+2(l-m-\epsilon))$, where $r\in\{1,\ldots,d-1\}$, $\epsilon\in\{0,1\}$ and $l\in\mathbb{N}$ with $l\geq\max(\epsilon,\frac{d+1}{2}-r+m)$.
The multiplicity of the eigenvalue corresponding to the choice of a triple $(r,\epsilon,l)$ as above is given by
\[\frac{d(\frac{d+1}{2}+r-m+2l-\epsilon)}{(r+l)(\frac{d+1}{2}-m+l-\epsilon)}\cdot\left(\begin{array}{c}d+l-\epsilon\\ d\end{array}\right)\cdot\left(\begin{array}{c}d-1\\ d-r-\epsilon\end{array}\right)\cdot\left(\begin{array}{c}\frac{d-1}{2}+r-m+l\\ d\end{array}\right).\]
\item $2l(2l+d-1-2m)$, where $l\in\mathbb{N}$, $l\geq\max(0,m+\frac{d+1}{2})$, with multiplicity
\[\prod_{k=2}^d(1+\frac{l}{k-1})\cdot(1+\frac{2l-\frac{d+1}{2}-m}{d})\cdot\prod_{j=2}^d(1+\frac{l-\frac{d+1}{2}-m}{d-j+1}).\]
\item $2(d+l)(d+1+2(l-m))$, where $l\in\mathbb{N}$, $l\geq\max(0,m-\frac{d+1}{2})$, with multiplicity
\[\prod_{k=2}^d(1+\frac{l}{k-1})\cdot(1+\frac{2l+\frac{d+1}{2}-m}{d})\cdot\prod_{j=2}^d(1+\frac{l+\frac{d+1}{2}-m}{d-j+1}).\]
\een
\end{ethm}

\noindent Note that, since $\CP^d$ is a symmetric space, the spectrum of every Dirac operator twisted with a homogeneous bundle over $\CP^d$ is symmetric about the origin. 
Hence the spectrum of the Dirac operator of $\CP^d$ twisted with $\gamma_d^m$ can be easily deduced from that of its square.\\

\noindent We point out that the computations done by M. Ben Halima in \cite[Thm. 1]{BenHalima08} contain a minor mistake (his $m$ should be replaced by $-m$).
It can be also checked that, up to a factor $4(d+1)$ (his convention for the Fubini-Study metric is different from ours), our values coincide with his (his $k$ is our $l$ and his $l$ is our $d-r$).\\

\noindent We can now formulate the 

\begin{ethm}\label{t:compspec}
Let $d<n$ be positive odd integers.
Then the spectrum of the square of the Dirac operator of $\CP^d$ twisted with the spinor bundle of the normal bundle of the canonical embedding $\CP^d\rightarrow\CP^n$ is given by the following family of eigenvalues:
\ben\item $2(r+l)\cdot(2d+1-n+2(s+l-\epsilon))$, where $r\in\{1,\ldots,d-1\}$, $s\in\{0,\ldots,n-d\}$, $\epsilon\in\{0,1\}$ and $l\in\mathbb{N}$ with $l\geq\max(\epsilon,\frac{n+1}{2}-r-s)$.
The multiplicity of the eigenvalue corresponding to the choice of a $4$-tuple $(r,s,\epsilon,l)$ as above is given by
\[\frac{d(d-\frac{n-1}{2}+r+s+2l-\epsilon)}{(r+l)(d-\frac{n-1}{2}+s+l-\epsilon)}\cdot\left(\begin{array}{c}n-d\\ s\end{array}\right)\cdot\left(\begin{array}{c}d+l-\epsilon\\ d\end{array}\right)\cdot\left(\begin{array}{c}d-1\\ d-r-\epsilon\end{array}\right)\cdot\left(\begin{array}{c}d-\frac{n+1}{2}+r+s+l\\ d\end{array}\right).\]
\item $4l(l+s+d-\frac{n+1}{2})$, where $s\in\{0,\ldots,n-d\}$, $l\in\mathbb{N}$, $l\geq\max(0,\frac{n+1}{2}-s)$, with multiplicity
\[\left(\begin{array}{c}n-d\\ s\end{array}\right)\cdot\prod_{k=2}^d(1+\frac{l}{k-1})\cdot(1+\frac{2l-\frac{n+1}{2}+s}{d})\cdot\prod_{j=2}^d(1+\frac{l-\frac{n+1}{2}+s}{d-j+1}).\]
\item $2(d+l)(2d-n+1+2(l+s))$, where $s\in\{0,\ldots,n-d\}$, $l\in\mathbb{N}$, $l\geq\max(0,\frac{n-1}{2}-d-s)$, with multiplicity
\[\left(\begin{array}{c}n-d\\ s\end{array}\right)\cdot\prod_{k=2}^d(1+\frac{l}{k-1})\cdot(1+\frac{2l+d-\frac{n-1}{2}+s}{d})\cdot\prod_{j=2}^d(1+\frac{l+d-\frac{n-1}{2}+s}{d-j+1}).\]
\een
\end{ethm}

\noindent{\it Proof}:
Recall that, by Corollary \ref{c:SigmaTperp}, there exists a unitary and parallel isomorphism
\[ \Sigma(T\CP^d)\otimes\Sigma(T^\perp\CP^d)\cong\bigoplus_{s=0}^{n-d}\left(\begin{array}{c}n-d\\ s\end{array}\right)\cdot\Sigma(T\CP^d)\otimes\gamma_d^{\frac{n-d}{2}-s},\]
where $\gamma_d$ is the tautological bundle of $\CP^d$ and $\left(\begin{array}{c}n-d\\ s\end{array}\right)$ stands for the multiplicity with which the subbundle $\Sigma(T\CP^d)\otimes\gamma_d^{\frac{n-d}{2}-s}$ appears in the splitting.
Therefore, the eigenvalues of the twisted Dirac operator acting on $\Sigma(T\CP^d)\otimes\Sigma(T^\perp\CP^d)$ are those of $\Sigma(T\CP^d)\otimes\gamma_d^{\frac{n-d}{2}-s}$, where $s$ runs from $0$ to $n-d$.
Moreover, the multiplicity of the eigenvalue corresponding to some $s$ is $\left(\begin{array}{c}n-d\\ s\end{array}\right)$ times the multiplicity computed in Lemma \ref{l:multvalp}.
Replacing $m$ by $\frac{n-d}{2}-s$, Theorem \ref{t:specDgammadm} gives the result.
\findemo
$ $\\

\noindent Note that $(d+1)^2$ is always an eigenvalue for the squared operator $(D_M^{\Sigma N})^2$: if $d=1$, take $s=\frac{n-1}{2}$ and $l=1$ in the second family of eigenvalues; if $d>1$, take $r=\frac{d+1}{2}$, $s=\frac{n-d}{2}$ and $\epsilon=0=l$ in the first family.\\

\noindent Using Theorem \ref{t:compspec}, we are now able to compute the smallest eigenvalue of the twisted Dirac operator:

\begin{prop} \label{pro:loweigen}
The lowest eigenvalue for the square of the Dirac operator of $\CP^d$ twisted with the spinor bundle of the normal bundle of the canonical embedding $\CP^d\rightarrow\CP^n$ is equal to $0$ for $d<\frac{n+1}{2}$ and to $(n+1)(2d+1-n)$ for $d\geq\frac{n+1}{2}$.
\end{prop}

\noindent{\it Proof}: Let us consider the first family of eigenvalues with $\epsilon=0$ (the same computation remains true for $\epsilon=1$). For $r+s\geq \frac{n+1}{2},$ which implies $d-\frac{n-1}{2}\leq r$, the minimum is attained for $l=0$ and we find the eigenvalues 
$2r(2d+1-n+2s),$ which are increasing functions with respect to $s$ with $s\geq \frac{n+1}{2}-r$. Here two cases occur: 
\begin{enumerate}
\item Case where $\frac{n+1}{2}-r\geq 0$, the eigenvalues become $4r(d+1-r)$ and we distinguish the two subcases: 
\begin{enumerate}
\item For $d\leq \frac{n+1}{2}$, then the lowest eigenvalue is equal to $4d.$
\item For $\frac{n+1}{2}<d$, the lowest eigenvalue is $(n+1)(2d+1-n)$.
\end{enumerate} 
\item Case where $\frac{n+1}{2}-r<0$ which implies $\frac{n+1}{2}< d$. Hence, the lowest eigenvalue is equal to $(n+1)(2d+1-n)$.
\end{enumerate}

\noindent Now for $r+s<\frac{n+1}{2}$, we take $l=\frac{n+1}{2}-r-s$. Thus the eigenvalues are equal  
$2(n+1-2s)(d+1-r)$ which are decreasing functions in $s$ with $0\leq s\leq \frac{n-1}{2}-r$.
We have:
\begin{enumerate}
\item Case where $\frac{n-1}{2}-r\leq n-d$. We then get the eigenvalues $4(1+r)(d+1-r)$ with $d-\frac{n+1}{2}\leq r\leq \frac{n-1}{2}$. Here two cases occur: 
\begin{enumerate}
\item For $d\leq \frac{n+1}{2}$, the lowest eigenvalue is equal to $8d$. 
\item For $d> \frac{n+1}{2}$, the lowest eigenvalue is equal to $(n+3)(2d+1-n).$
\end{enumerate}
\item Case where $\frac{n-1}{2}-r> n-d$, we get the eigenvalues $2(2d-n+1)(d+1-r)$ with $1\leq r\leq d-\frac{n+3}{2}$. In this case, we have that $d>\frac{n+1}{2}$ and the lowest eigenvalue is equal to $(n+5)(2d-n+1)$.
\end{enumerate}
For the second family of eigenvalues, we distinguish the cases: 
\begin{enumerate}
\item Case where $\frac{n+1}{2}-s\leq 0$ which implies that $d\leq \frac{n-1}{2}$, we take $l=0$. The lowest eigenvalue is then equal to $0$.
\item Case where $\frac{n+1}{2}-s>0$. The eigenvalues become $2d(n+1-2s)$ with $0\leq s\leq \frac{n-1}{2}$. Two cases occur 
\begin{enumerate}
\item For $d\leq \frac{n+1}{2}$, the lowest eigenvalue is $4d.$
\item For $d>\frac{n+1}{2}$, the lowest eigenvalue is $2d(2d+1-n)$
\end{enumerate}
\end{enumerate}
For the last family of eigenvalues, we consider the two cases: 
\begin{enumerate}
\item Case where $\frac{n-1}{2}-d-s>0$, which implies that $d<\frac{n-1}{2}$, we take $l=\frac{n-1}{2}-d-s$. We find the lowest eigenvalue $0$ after substituting. 
\item Case where $\frac{n-1}{2}-d-s\leq 0$. In this case $l=0$ and we get $2d(2d-n+1+2s)$. Here two cases occur: 
\begin{enumerate} 
\item For $d> \frac{n-1}{2}$, the lowest eigenvalue is $2d(2d-n+1)$. 
\item For $d\leq\frac{n-1}{2}$, the lowest eigenvalue is $0$. 
\end{enumerate}
\end{enumerate}
\hfill$\square$\\

\noindent Next we show that the estimate (\ref{eq:majalphareel}) is not always sharp. 
We consider the simplest case where $d=1$ and compare the mul\-ti\-pli\-ci\-ties of the eigenvalues $0$ and $4$ with $2\left(\begin{array}{c}n\\ \frac{n+1}{2}\end{array}\right)$, which is the \emph{a priori} number of eigenvalues bounded by $4$ in (\ref{eq:majalphareel}).
The multiplicity of the eigenvalue $0$ is equal to 
$$\sum_{s=0}^{\frac{n-3}{2}} \left(\begin{array}{c}n-1\\ s\end{array}\right)(\frac{n-1}{2}-s)+\sum_{s=\frac{n+1}{2}}^{n-1} \left(\begin{array}{c}n-1\\ s\end{array}\right)(s-\frac{n-1}{2}),$$
which is equal to $\sum_{s=0}^{\frac{n-3}{2}} \left(\begin{array}{c}n-1\\ s\end{array}\right)(n-1-2s)$ since by replacing $s$ by $(n-1)-s$ the second sum  is equal to the first one.
A short computation gives $\sum_{s=0}^{\frac{n-3}{2}} \left(\begin{array}{c}n-1\\ s\end{array}\right)(n-1-2s)=\frac{n-1}{2}\cdot\left(\begin{array}{c}n-1\\ \frac{n-1}{2}\end{array}\right)$.
On the other hand, the multiplicity of the eigenvalue $4$ is equal to $4\left(\begin{array}{c}n\\ \frac{n-1}{2}\end{array}\right)$.
Hence the sum of these two multiplicities is $(\frac{n-1}{2}+4)\cdot\left(\begin{array}{c}n-1\\\frac{n-1}{2}\end{array}\right)$. 
That number is always greater than $2\left(\begin{array}{c}n\\
\frac{n+1}{2}\end{array}\right)$.
However, if the multiplicity of the eigenvalue $0$ is smaller than $2\left(\begin{array}{c}n\\ \frac{n+1}{2}\end{array}\right)$ for $n=3,5,7$, it is greater for $n\geq 9$.
Thus, the equality in (\ref{eq:majalphareel}) is optimal for $n=3,5,7$ but is \emph{never} optimal as soon as $n\geq 9.$
In particular, the twisted Dirac operator on K\"ahler submanifolds behaves very differently from that on submanifolds immersed in real spaceforms, where analogous upper bounds are sharp in any dimension.

\providecommand{\bysame}{\leavevmode\hbox to3em{\hrulefill}\thinspace}

\end{document}